\documentclass[a4paper,10pt]{article}

\def \Z {{\mathbf {Z}}}
\def \R {{\mathbf {R}}}

\def\u{\bigsqcup}
\def\eps{\varepsilon}
\title{ On Mixing Rank One Infinite Transformations }
\author{V.V. Ryzhikov\footnote{\large This work is partially  supported by  the grant
NSh 8508.2010.1.}}

\begin{document}
\Large
\maketitle

D.Ornstein has defined rank one transformations; he has shown that they
could be mixing and that they then only commuted with their powers \cite{O}.
J.L. King has later proved that mixing rank one transformations had
the much stronger property of being MSJ
 (see \cite{K1}).  Examples of infinite transformations with MSJ  were obtained by J. Aaronson and  M. Nadkarni in \cite{AN}. 
J.-P. Thouvenot and the author showed via different approaches that the centralizer of a mixing rank-one {\it infinite} measure preserving transformation was trivial. Thouvenot used Ornstein's method.  In this note the author presents his joining  proof based on the technique of \cite{R}. We also consider  constructions with algebraic spacers as well as   a class of
"Sidon constructions"  to produce new examples of mixing rank one transformations.  In connection with 
Gordin's question on the existence of homoclinic ergodic actions for a zero entropy system \cite {G},\cite{K2} we also discuss  Poisson suspensions of some modifications of Sidon rank one  constructions.

\section{ Mixing rank one constructions and some applications.} 
We recall the notion of rank one transformation. Let us 
consider an infinite (or finite) Lebesgue space $(X,\mu)$.
An automorphism (a measure-preserving invertible transformation)
$T:X\to X$ is said to be of
{\it rank one}, if there is a sequence $\xi_j$ of measurable partitions
of $X$ in the form
$$\xi_j= \{ E_j,\ TE_j, \ T^2 E_j,\ \dots, T^{h_j}E_j, \tilde{E}_j\}.$$
such that $\xi_j$ converges to the partition onto points.
The collection $$E_j, TE_j, \ T^2 E_j,\ \dots, T^{h_j}E_j$$ is called
Rokhlin's tower. We put 
$$U_j=E_j\u TE_j\u T^2 E_j\u\dots\u T^{h_j}E_j$$ 
( $\tilde{E}_j$ is the set
$X\setminus \u_{i=0}^{h_j} T^i E_j$).

{\bf Mixing.}
An infinite transformation $T$ is called mixing if for any $\eps >0$ and  sets $A,B$ of finite measures there is $N$ such that for all $|k|>N$ 
$$
\mu(T^kA\cap B) < \eps .$$

To construct  examples of mixing rank one infinite transformations one can use the methods by D. Ornstein or T. Adams (see \cite{A},\cite{DR}).
Below we consider  more simple  examples.  

{\bf Rank one construction} is detrmined by $h_1$ and a sequence $r_j$  of cuttings and  a sequence $\bar s_j$ of spacers
$$ \bar s_j=(s_j(1), s_j(2),\dots, s_j(r_j-1),s_j(r_j)).$$  We recall its definition.
Let our  $T$ on the step $j$  be assosiated with  a collection of disjoint sets (intervals)
$$E_j, TE_j T^2,E_j,\dots, T^{h_j}E_j.$$
We cut $E_j$ into $r_j$ sets (subintervals)  of the same measure
$$E_j=E_j^1\u E_j^2\u  E_j^3\u\dots\u E_j^{r_j},$$  
then for all $i=1,2,\dots, r_j$ we  consider columns
$$E_j^i, TE_j^i ,T^2 E_j^i,\dots, T^{h_j}E_j^i.$$
Adding $s_j(i)$ spacers we obtain  
$$E_j^i, TE_j^i T^2 E_j^i,\dots, T^{h_j}E_j^i, T^{h_j+1}E_j^i, T^{h_j+2}E_j^i, \dots, T^{h_j+s_j(i)}E_j^i$$
(the above intervals   are disjoint).
For all  $i<r_j$ we set
$$TT^{h_j+s_j(i)}E_j^i = E_j^{i+1}.$$ 
Now we obtaine a tower 
$$E_{j+1}, TE_{j+1} T^2 E_{j+1},\dots, T^{h_{j+1}}E_{j+1},$$
where 
 $$E_{j+1}= E^1_j$$
$$T^{h_{j+1}}E_{j+1}=T^{h_j+s_j(r_j)}E_j^{r_j} $$
$$h_{j+1}+1=(h_j+1)r_j +\sum_{i=1}^{r_j}s_j(i).$$
So we have described a general rank one construction.

{\bf Algebraic spacers (a.s.) constructions.}  Let $r_j$ be prime, $r_j\to\infty$.  We fix and consider  generators $q_j$ in the multiplicative groups (associated with the sets $\{1, 2,\dots, r_j-1\}$) of the  fields $\Z_{r_j}$.
For some sequence  $\{H_j\}$, $H_j\geq r_j$, we define a spacer sequence
$$s_j(i)=H_j + \{q_j^i\} -\{q_j^{i+1}\},  \ \ i=1, 2, \dots,  r_j-1,$$
where $H_j$  provides  the infinity of our  measure space,  $s(r_j)$ is arbitrary,
$\{q^i\}$ denotes the number from $\{1, 2,\dots, r_j-1\}$ that corresponds to $q^i$.
This rank one a.s. construction will be mixing.
To prove the mixing we must check and apply two properties of our spacers: for  
$n<r_j$ we have 

(1) $ - r_j \leq S_j(i,n):=\sum_{k=1}^{n} s_j(i+k)  - nH_j\leq r_j , \ \  i=1, 2, \dots, r_j-n-1$  (Ornstein's property);

(2)  for $ i\in\{1, 2, \dots, r_j-n-1\}$ all values $ S_j(i,n)$ are different (a nice injectivity property).

Indeed,  $ S_j(i,n)=  S_j(m,n)$  implies $q^i-q^{i+n}=q^m-q^{m+n}$, $q^i=q^m$, $i=m$.

 For any  sequence $n_j<(1-\eps)r_j$ ($\eps> 0$ is fixed) and any  positive $f\in L_2$ we get 
$$\Vert P(j, n_j)f\Vert_{L_2} \to 0,   \ j\to \infty,$$
where
$$  P(j,n_j):=\frac{1}{r_j-n_j-1}\sum_{i=1}^{r_j-n_j-1} T^{ S_j(i,n_j)}.$$
Indeed,  from (1) and (2) we get 
$$ \Vert P(j, n_j)f\Vert_{L_2}\leq \frac{1}{\eps r_j} \Vert\sum_{s=-r_j}^{r_j} T^{ s} f\Vert_{L_2}\ \ \to \ 0.$$
Let $$h_j \leq m_j=n_j(H_j+h_j) +t_j<h_{j+1},  \ \ 0\leq t_j\leq H_j+h_j,$$ 
then standard rank one
method of estimations of mixing via averaging (see \cite{A},\cite{DR}) gives
$$\langle T^{m_j} \chi_A, \ \chi_B\rangle \leq   \ (\|  P(j,n_j)\chi_A\| + \|P(j,n_j+1)\chi_A\|\ +\ \| P(j+1,1)\chi_A\| + \eps) \|\chi_B\| $$
for arbitrary $\eps>0$ and all large $j$. Thus, for $m\to\infty$ we get 
$$\langle T^{m} \chi_A, \ \chi_B\rangle \to 0.$$

{\bf REMARK.} In fact we can use a.s. constructions  (as $H_j$ and $s_j(r_j)$ are sufficiently small) for  examples of mixing rank one transformations of a   Probability  space as well.  We follow Ornstein's method but stochastic spacers are replaced 
now by algebraic ones.  

Below we consider  mixing constructions only  for infinite measure spaces but with an obvious proof of the mixing property.

{\bf Sidon constructions.} Let $r_j\to\infty$ and for each step $j$
$$ h_j<<s_j(1)<< s_j(2)<<\dots<< s_j(r_j-1)<<s_j(r_j) \eqno (\ast).$$ 
Then for fixed $\xi_{j_0}$-measurable $A,B\in U_{j_0}$   we have 
$$\mu(A\cap T^mB)\leq \mu(A)/r_j$$
as $m\in [h_j, h_{j+1}],$ $j>j_0$.  
Thus, for any $A,B$ of finite measure we get
$$\mu(A\cap T^mB)\to 0.$$

  A reader can check that our construction has "Sidon property":  only one column in $U_j$  may contain an intersection
$U_j\cap T^mU_j$ as $h_{j+1}>m>h_j$. 
(We recall that Sidon set is a subset $S$ of $\{1,2,\dots,N\}$ such that  $S\cap S+m$ ($m>0$) may contain not more than one point.
Its cardinality can be a bit grater than $\sqrt{N}$.) Any Sidon construction (not necessarily satisfied  ($\ast$)) have to be  mixing. 
    There are Sidon constructions (with $r_j >> h_j$ but with "minimal"  spacers) with a  rate of correlations 
 $$\mu(A\cap T^mA)\leq C \frac{\psi(m)}{\sqrt{m}} ,$$
where $\psi(m)$  satisfies   the condition
$$\frac{\psi( h_{j+1})}{\sqrt{ h_{j+1}}  } \leq   \frac{\psi(m)}{\sqrt{m}}, \ \ m\in [h_j+1, h_{j+1}],$$
and $\psi(m)\to \infty$ as slowly as we want (for example,   $\psi(m)=\ln \ln{(m)}$). 
 
Indeed, given $h_j$   we choose   $h_{j+1}\sim  h_j r_j^2$  such that  $\psi(h_{j+1})\geq \sqrt{h_j}$,
then for all $m\in [h_j+1, h_{j+1}] $ we have 
$$\mu(A\cap T^mA)\leq \mu(A)/r_j \leq  C  \frac{\sqrt{h_j}}{\sqrt{ h_{j+1}}}= C
\frac{\psi( h_{j+1})}{\sqrt{ h_{j+1}}  }\frac{\sqrt{h_j}}{\psi({ h_{j+1}})}\leq C \frac{\psi(m)}{\sqrt{m}}.$$

REMARK. In \cite{P} A. Prikhod'ko  proposed non-trivial constructions of automorphisms of a Probability space with simple spectrum and the  correlation sequence 
$$\langle T^mf, f\rangle = O(m^{-1/2+\eps}).$$  We obtained  simply the same result
for some Sidon constructions  due to  a spacer freedom in the infinite measure case.

Adding to our Sidon  constructions  a special (vanishing non-mixing) parts as in \cite{DR}
we provide  simple spectrum for the operator $exp(T)$ and 
$$\mu(A\cap T^mB)\leq  \eps_j +\mu(A)/r_j,$$
where $\eps_j\to 0$ very very slowly.
Following this way we get again one of  the  results of \cite{DR}: there are mixing Poisson suspensions  of simple spectrum. 

{\bf  Homoclinic groups.} The above transformations $T$ can possess  ergodic homoclinic groups $H(T)$.  Following M.Gordin (see \cite{G},\cite{K2})  we recall that 
 $$H(T)=\{ S\in Aut(X,\mu): \ T^nST^{-n}\to Id,  n\to\infty\}.$$
We see that   all $S$ with $\mu(supp(S))<\infty$ are in $H(T)$.  
This observation suggest us a construction of  a dissipative transformation $S\in H(T)$. 
Indeed, for $T$ we can build an almost transversal transformation $S$ which is more and more close  to the identity on $j$-spacers
as $j\to \infty$.  If we fix a finite measure set $A\subset U_j$, then $\mu(A\Delta T^nST^{-n}A)\to 0$ since $T^nST^{-n}$
will be close to the identity on $A$. It is not so hard to provide   $S$ to be dissipative.
( for  $n>h_j$  we have $\mu(A\Delta T^nST^{-n}A)< 2\eps_j +2\mu(A)/r_j +\delta_j,$ where $\eps_j$ as above, and some $\delta_j\to 0$ because our $S$  is more and more close  to the identity on $j$-spacers).

Thus  we obtain  the Poisson suspension $T_\ast$ with simple spectrum possessing the homoclinic Bernoulli transformation $S_\ast$
(this is somewhat in contrast with the natural situation in which  Bernoulli transformations play the role of $T$). 
It gives  a new  solution of Gordin's  question  with respect to King's one\cite{K2}. 
We note also that {\it a mixing rank one transformation of a  finite measure space
has no homoclinic transformations $S\neq Id$}. 

M. Gordin also asked the author on an example of a zero entropy transformation
with a homoclinic ergodic flow. 
 Now we give a  short  solution based on Poisson
suspensions.
Let $\tilde T$ a mixing  rank one transformation of infinite measure space
$( \R^+,  \mu)$.
We set $X=\R\times  \R^+$    and define $T:X\to X$ by the equality
$$T(x,y)=(x,\tilde T (y)).$$   It possesses a  homoclinic flow.  Indeed,
let 
$$  S_t(x,y)=(x+\varphi(y)t, y),$$
where $\varphi >0$ and $\varphi(y)\to 0$ as $y\to\infty$.
The flow $S$ is dissipative and homoclinic: for  $f\in L_2( \R\times  \R^+,\mu\times\mu)$ we have 
$$\|f - T^nS_tT^{-n} f\|_2\to 0,\ \ n\to\infty. $$

The latter is obvious for $f=\chi_{I_1\times I_2}$ ($I_1, I_2$ are finite intervals), so it is true for all $f\in L_2$.  
The flow $S_{\ast t}$ is  Bernoulli flow with infinite entropy.    

The Poisson suspensions $T_\ast$ can be of zero entropy.  For this we use $\tilde T$ with $exp(\tilde T)$ of simple spectrum.  The spectral type of $T_\ast$ is singular, it coincides with the spectrum of $exp(\tilde T)$, but  the  spectrum of $T_\ast$ has  infinite multiplicity. 

{\bf REMARK.}  In connection with the work \cite{FL} let us note that the Poisson suspensions
could be a nice source   of self-similar flows (with Lebesgue spectrum). 
Let $\tilde T_t$ a  flow of infinite measure space
$( Y,  \mu)$.
We set again $X=\R\times  Y$    and define $T_t:X\to X$ by the equality
$$T_t(x,y)=(x,\tilde T_{e^xt} (y)).$$  
Let us consider a flow $G_a$
$$  G_a(x,y)=(x+a, y).$$
We get 
$$G_aT_tG_{-a}=T_{e^{-a}t},$$
$$G_{\ast a}T_{\ast t}G_{\ast -a}=T_{\ast e^{-a}t}.$$
We get now that   all $S$ with $\mu\times\mu(supp(S))<\infty$ are in $H(T_t)$.  
Thus our flows $T_t$ and $T_{\ast t}$ possess ergodic homoclinic groups.
The author is indebted to A.I.Danilenko and E.Roy for arguments (based on \cite{J-R}) showing  
that a large class of flows $T_{\ast t}$ as above could be of zero entropy.

M. Gordin proved that a transformation with an ergodic homoclinic group have to be 
mixing.  We note a bit more: 

{\bf Assertion.} \it An  action with  an ergodic homoclinic group
is mixing of all orders.  \rm 

Thus,  Ledrappier's non-multiple-mixing actions (with Lebesgue spectrum) have no  ergodic homoclinic groups. 

The assertion follows from 

{\bf Lemma.} \it  Let $T_i$ and $T_j$ be sequences of elements of a measure-preserving action on a Probability space, let $S$ satisfy 

(1) $T^{-1}_iST_i\to Id$, $T^{-1}_jST_j\to Id$.

 If 
for  a measure $\nu$   on $X\times X\times X$ 
and  measurable sets  $A,B,C$  it holds
$$\mu (A\cap T_i B\cap T_j C) \to \nu (A\times B\times C),$$ then   
$$\nu (SA\times B\times C)=\nu (A\times B\times C).$$

If in addition $\nu(X\times B\times C)=\mu (B)\mu(C)$ for all $A,B$
and all $S$  satisfied (1) generate together an ergodic action, then
$\nu=\mu\times\mu\times\mu$.
\rm

Proof is obvious:
$$\nu(A\times B\times C)=\lim_{i,j}\mu (A\cap T_i B\cap T_j C)=\lim_{i,j}\mu (SA\cap ST_i B\cap ST_j C)=$$
$$=\lim_{i,j} \mu (SA\cap T_iT^{-1}_iST_i B\cap T_jT^{-1}_jST_j C)= 
\nu (SA\times B\times C).$$

If $\nu(X\times B\times C)=\mu (B)\mu(C)$ for all $A,B$, then 
$\nu$ is a joining of an ergodic action with the identity action on $X\times X$.
Thus, $\nu$ have to be the direct  product of the projections  $\mu$ and $\mu\times\mu$.

\section{Self-joinings of mixing rank one transformations}
A {\it self-joining} (of order 2) is a $T\times T$-invariant
measure $\nu$ on $X\times X$ with marginals equal to $\mu$:
$$\nu(A\times X)=\nu(X\times A)=\mu(A).$$
Off-diagonal measures $\Delta^k$ are defined by the equality
$$\Delta^k(A\times B)=\mu(T^kA\cap B).$$

We say that $T$ has {\it minimal self-joinings} ({\it MSJ})
if all its ergodic self-joinings are off-diagonals $\Delta^k$ (a joining is called ergodic if it is not a convex sum of two different joinings).
\vspace{3mm}

\bf THEOREM 1. \it Mixing rank-one infinite transformations have MSJ.

\bf COROLLARY. \it
A mixing rank-one infinite transformation commutes only with its powers. 
\rm
\vspace{3mm}

{\bf REMARK.} It follows from corollary and the results of  \cite{Ro} that
the Poisson suspensions of the mixing rank one transformations constructed in  \cite{DR} have trivial centralizer. What is proved here is a partial answer to the question of E. Roy as
to whether the Weak Closure Theorem  of J. King extends to infinite
measure rank one transformations.
Theorem 1  in the finite measure case was proved by  J.~King (see \cite{K1},\cite{R}). 
  
\vspace{3mm}

Proof of Corollary. Suppose that an automorphism $S$ commutes with $T$.
The joining $\Delta_S=(Id\times S)\Delta$ is ergodic:
the system $(T\times T, X\times X, \Delta_S)$ is isomorphic to
the ergodic system $(T, X,\mu)$. Then for some $i$ we get
$(Id\times S)\Delta =(Id\times T^i)\Delta,$ so $S=T^i$.
\vspace{3mm}

\bf LEMMA. \it Let $\nu$ be a joining of a rank one transformation $T$. Suppose that $\nu$ is disjoint from all $\Delta^k$.
There is a sequence $a^k_j$ such that for all $j$ we have $\sum_k a^k_j \leq 2$ and the inequality
$$\nu(A\times B) \leq \limsup_j \sum_k a^k_j \Delta^k(A\times B)$$
holds for all finite measure sets $A,B$. 
\rm
\vspace{3mm}

Proof of Lemma.

We define  sets $C^k_j$ and  measures $\Delta^k_j$
:
$$ C^k_j=\u_{i=0}^{h_j-k} T^i E_j\times T^{i+k} E_j, \ \ 
k\in[0,h_j],$$

$$ C^k_j=\u_{i=0}^{h_j+k}T^{i-k} E_j\times T^i E_j, \ \  k\in[-1,-h_j],
$$
$$\Delta^k_j (A\times B)= \Delta^k_j ((A\times B)\cap C^k_j),\ \  k\in[-h_j,h_j].$$
 For  $k\in [0,h_j]$ we denote by $N(k, A,B)$  the number of $i$-s r such that the corresponding  blocks 
$T^i E_j\times T^{i+k}E_j$ are in $A\times B$.
 For $k\in[-1,-h_j]$ we  consider the blocks  $T^{i -k}E_j\times T^{i}E_j$. Thus,  $N(k, A,B)$ is the number of
blocks in $C^k_j\cap (A\times B)$.

Let us define for $k\in [0, h_j]$
$$ a_j^k =\frac { \nu( E_j\times T^k E_j)}{\mu(E_j)}$$
and for $k\in [- h_j, -1]$
$$ a_j^k =\frac { \nu( T^{-k}E_j\times E_j)}{\mu(E_j)}.$$

Using the  invariance $\nu( T^{-k+m}E_j\times T^m E_j)=
\nu( T^{-k}E_j\times E_j)$ we get
$$\nu(A\times B)=\mu(E_j)\sum_{|k|\leq h_j} a_j^kN(k, A,B) =
\sum_{|k|\leq h_j} a_j^k\Delta^k_j(A\times B).$$
From this we obtain
$$
\sum_{|k|\leq h_j} a_j^k\Delta^k_j \to \nu$$
since arbitrary measurable sets $A', B'$
can be approximated by $\xi_j$-measurable ones. 
We get
$$\nu(A\times B) \leq \lim_j \sum_k a^k_j \Delta^k_j(A\times B)\leq \limsup_j \sum_k a^k_j \Delta^k(A\times B)$$
and 
$$\sum_{|k|\leq h_j} a_j^k\leq 2.$$ 
The latter follows from marginal properties of a joining:  
$$\mu(E_j)\sum_{k=0}^{h_j} a^k_j=\sum_{k=0}^{k=h_j}\nu( E_j\times T^k E_j)\leq \nu(E_j\times X) ={\mu(E_j)}.$$ 
Lemma is proved.

To prove Theorem we apply Lemma 
and a consequence of the mixing: for any $\eps>0$ there is $M$
such that for all  $j$
$$
\sum_{|k|>M}a^k_j\Delta^k(A\times B)=
\sum_{|k|>M}a^k_j\mu(T^kA\cap B)<2\eps.$$

If $\nu(A\times B)> 0$, then $$\limsup_j \sum_k a^k_j \Delta^k(A\times B)>0,$$ and we find  $k'$ such that
$a^{k'}_j$ is not vanishing as $j\to\infty$, so
$$\nu \geq a\Delta^{k'}, \ \  a>0,$$ hence, $\nu=\Delta^{k'}$.

\vspace{5mm}

The author is thankful to M. Gordin and J.-P. Thouvenot  for stimulating questions
and discussions.

 E-mail: vryz@mail.ru
\end{document}